\documentclass[12pt]{amsart}

\usepackage{amssymb}
\usepackage{latexsym}
\usepackage{amsmath}

\def\kh{\mathbb{H}}
\def\dh{\widehat{\kh}}

\def\kl{{\mathbb L}}

\def\c{{\mathbb C}}

\def\l{{\mathcal L}}
\def\b{{\mathcal B}}

\def\sp{\overline{sp}^w}

\def\alg{{\mathcal Alg}}

\begin{document}

\newtheorem{defi}{Definition}[section]
\newtheorem{prop}{Proposition}[section]
\newtheorem{theo}{Theorem}[section]
\newtheorem{lemm}{Lemma}[section]
\newtheorem{coro}{Corollary}[section]

\title
[Compact Kac algebras and commuting squares]
{Compact Kac algebras and commuting squares}
\author{Teodor Banica}
\address{
C.N.R.S., Institut de Math\'ematiques de Jussieu, Paris and\newline 
Departement of Mathematics, University of California, Berkeley, CA 94720.}
\email{banica@math.jussieu.fr, banica@math.berkeley.edu}

\begin{abstract}
We consider commuting squares of finite dimensional von Neumann algebras having the algebra of complex numbers in the lower left corner. Examples include the vertex models, the spin models (in the sense of subfactor theory) and the commuting squares associated to finite dimensional Kac algebras. To any such commuting square we associate a compact Kac algebra and we compute the corresponding subfactor and its standard invariant in terms of it.
\end{abstract}

\maketitle

\section{Introduction}

There exist several constructions of subfactors using quantum groups and vice versa. The precise relation between this objects is far from being clear. In this paper we give an answer for certain subfactors and quantum groups, related to vertex models.

A vertex model (in the sense of subfactor theory) is a commuting square of the form
$$\begin{matrix}
{\c}\otimes M_k(\c ) &\subset & M_n(\c )\otimes M_k(\c )\cr
\cup &\ &\cup\cr
{\c} &\subset &u(M_n(\c )\otimes{\c})u^*\end{matrix}$$
where $u\in M_n(\c )\otimes M_k(\c )$ satisfies the biunitarity condition: $u$ and $(t\otimes id)u$ are
unitaries, where $t$ is the transposition of $M_n(\c )$.

A spin model (in the sense of subfactor theory) is a commuting square of the form
$$\begin{matrix}
\Delta &\subset & M_n(\c )\cr
\cup &\ &\cup\cr
\c &\subset &w\Delta w^*\end{matrix}$$
where $\Delta\subset M_n(\c )$ is the algebra of diagonal matrices and
where $w\in M_n(\c )$ is a complex Hadamard matrix, i.e. a unitary all
whose entries are of modulus $n^{-1/2}$.

See e.g. \cite{js}, \cite{j3} for vertex models and spin models.

Recall also that associated to any finite quantum group (= dual of a finite dimensional Kac algebra, in this paper) $G$ is the following commuting square.
$$\begin{matrix}
C(G) &\subset & \l (l^2(G))\cr
\cup &\ &\cup\cr
\c &\subset & C(\widehat{G})\end{matrix}$$

These are all examples of commuting squares having $\c$ in the lower left corner.

The main result in this paper is a structure result, in terms of a compact quantum group of Kac type, for the subfactor associated to such a commuting square.

\begin{theo}
The vertical subfactor associated to a commuting square of finite dimensional von Neumann algebras of the form
$$\begin{matrix}
A &\subset & X\cr
\cup &\ &\cup\cr
\c &\subset & B\end{matrix}$$
is of the form $R\subset (B\otimes (R\rtimes\widehat{G}))^G$, where $G$ is a compact quantum group of Kac type, $R$ is the hyperfinite $II_1$ factor, $G\to Aut(B)$ is an action which is ergodic on the center of $B$ and $\widehat{G}\to Aut(R)$ is an outer action.
\end{theo}

The subfactor $R\subset (B\otimes (R\rtimes\widehat{G}))^G$ is a fixed point subfactor in the sense of \cite{kac}. By \cite{kac} its standard invariant is the Popa system of the action of $G$ on $B$. For vertex models this generalises our computation in \cite{ver}.

The quantum group $G$ is the one associated as in \cite{ver} to the vertex model obtained by performing two basic constructions to the commuting square. Note that for spin models we get compact quantum groups acting on $n$ points. See Wang \cite{wang2} and Bichon \cite{bi} for other examples of such exotic quantum groups.

The outer action of the discrete quantum group $\widehat{G}$ on $R$ is by construction a product type action. Since any finite quantum group arises from a vertex model, this generalises Yamanouchi's construction \cite{y}.

An important step in the proof is the fact that the commuting square in the theorem is isomorphic to a commuting square of the form
$$\begin{matrix}
A &\subset &  (B\otimes (A\rtimes\widehat{G} ))^G\cr
\cup &\ &\cup\cr
\c &\subset & B\end{matrix}$$
where $\widehat{G}\to Aut(A)$ is a certain action. This could be interpreted as giving a structure result for the commuting squares having $\c$ in the lower left corner.

More generally, we construct in \S 2 below commuting squares of the form
$$\begin{matrix}
(B_0\otimes (A_1\rtimes\widehat{G} ))^G &\subset &  (B_1\otimes (A_1\rtimes\widehat{G} ))^G\cr
\cup &\ &\cup\cr
(B_0\otimes (A_0\rtimes\widehat{G} ))^G &\subset & (B_1\otimes (A_0\rtimes\widehat{G} ))^G\end{matrix}$$

Such a commuting square has the property that each of its four algebras is endowed with its canonical trace in the sense of \cite{aut}. We do not know if the converse holds.

The rest of the paper is organised as follows. In \S 2 we associate commuting squares to compact quantum groups of Kac type. In the first appendix \S 4 we discuss the duality for vertex models. In the second appendix \S 5 we prove the result on the outerness of product type actions. In \S 3 we prove the main theorem, by using \S 2, \S 4, \S 5.

Part of this work was done during my participation at the Operator Algebra Program at the Erwin Schr\" odinger Institute for Mathematical Physics at Vienna and I would like to thank the Institute for financial support and J. Cooper and F. Lehner for the invitation. Another part of it was done during my visit at the University of California at Berkeley and I would like to thank D. Voiculescu for the invitation.

\section{Commuting squares of fixed point algebras}
 
Let $\kh$ be a compact Kac algebra with comultiplication $\Delta$ and antipode $S$. 
Denote by $\kh_\sigma$ the Kac algebra $(\kh
,\sigma\Delta ,S)$, where $\sigma$ is the flip. If $\beta :B\to
B\otimes\kh$ is a coaction on a finite dimensional finite von
Neumann algebra and $\pi :P\to
P\otimes\kh_\sigma$ is a coaction on a finite
von Neumann algebra define a
linear map
$$\beta\odot\pi :B\otimes P\to B\otimes P\otimes\kh_\sigma ,\,\,\, b\otimes
p\mapsto \pi (p)_{23}((id\otimes S)\beta (b))_{13}$$
Then its fixed point linear space
$$(B\otimes P)^{\beta\odot\pi}=\{ x\in B\otimes P\mid
(\beta\odot\pi )x=x\otimes 1\}$$
is a von Neumann subalgebra of $B\otimes P$ (cf. theorem 3.1 in \cite{kac}). The map $\beta\odot\pi$ is coassociative with respect to the comultiplication $\sigma\Delta$ of $\kh_\sigma$, but is not multiplicative in general; however, its fixed point algebra $(B\otimes P)^{\beta\odot\pi}$ has good properties. See \cite{kac} for more comments and for examples.

In this paper we consider only the case where $\pi$ is a dual coaction. Note that proposition 2.1 below provides an equivalent definition of $(B\otimes P)^{\beta\odot\pi}$ in this case.

Let $\int :\kh\to \c$ be the Haar functional, let $l^2(\kh )$ be its $l^2$-space and let $\dh\subset \b (l^2(\kh
 ))$ be the dual Kac algebra. If $\alpha :A\to A\otimes\dh$ is a coaction of $\dh$ on a finite von Neumann
  algebra $A$, the crossed
  product $A\rtimes_\alpha\dh$ is the von Neumann
  subalgebra of $A\otimes \b (l^2(\kh ))$ generated by $\alpha (A)$ and by
$1\otimes\kh$. There exists a unique coaction $\widehat{\alpha}$ of $\kh$
  on $A\rtimes_\alpha\dh$ such that
  $(A\rtimes_\alpha\dh )^{\widehat{\alpha}} =\alpha (A)$ and such that the
  copy $1\otimes\kh$ of $\kh$ is equivariant.
  
\begin{prop}
Let $\kh$ be a compact Kac algebra. If $\beta :B\to B\otimes\kh$ is a coaction on a
finite dimensional finite von Neumann algebra and $\alpha
:A\to A\otimes\dh_\sigma$ is a coaction on a finite von Neumann
algebra then
$$(B\otimes (A\rtimes_\alpha\dh_\sigma ))^{\beta\odot\widehat{\alpha}}=\sp \{ \beta (B)_{13}\cdot \alpha (A)_{23}\}$$
as linear subspaces of $B\otimes A\otimes {\b}(l^2(\kh_\sigma ))$. Moreover, the following diagram
$$\begin{matrix}
 \alpha (A)_{23} &\subset & (B\otimes (A\rtimes_{\alpha}\dh_\sigma 
))^{\beta\odot\widehat{\alpha}}\cr 
\cup &\ &\cup\cr 
{\c} &\subset & \beta (B)_{13}
\end{matrix}$$
is a non-degenerate commuting square of finite von Neumann algebras.
\end{prop}

\begin{proof}
By definition of the crossed product $A\rtimes_\alpha\dh_\sigma$ we have the following
equalities between subalgebras of $B\otimes A\otimes {\b}(l^2(\kh_\sigma
))$:
$$B\otimes (A\rtimes_\alpha\dh_\sigma )=B\otimes (\sp\{
\alpha (A)\cdot (1\otimes\kh_\sigma )\} ) =\sp\{ (B\otimes\kh_\sigma )_{13}\cdot
\alpha (A)_{23}\}$$
On the other hand since the coactions on finite dimensional algebras
are non-degenerate we have $B\otimes\kh_\sigma =sp\{ (1\otimes\kh_\sigma )\cdot
\beta (B)\}$ (see e.g. lemma 5.2 in \cite{kac}). Thus
$$B\otimes (A\rtimes_\alpha\kh_\sigma )= \sp\{ (1\otimes
1\otimes\kh_\sigma )\cdot \beta (B)_{13}\cdot \alpha (A)_{23}\}$$
Let us compute the restriction of $\beta\odot\widehat{\alpha}$ to $1\otimes
1\otimes\kh_\sigma$, to $\beta (B)_{13}$ and to $\alpha (A)_{23}$.

(a) the restriction of $\beta\odot\widehat{\alpha}$ to $1\otimes 1\otimes\kh_\sigma$ is $1\otimes
1\otimes\sigma\Delta$. In particular $\beta\odot\widehat{\alpha}$ has no fixed
points in $1\otimes 1\otimes\kh_\sigma$.

(b) the algebra $\alpha (A)_{23}$ is fixed by
$\beta\odot\widehat{\alpha}$.

(c) we prove that the algebra $\beta (B)_{13}$ is also fixed by
$\beta\odot\widehat{\alpha}$. For, let $\{ u_{ij}\}$ be an
orthonormal basis of $l^2(\kh_\sigma)$ consisting of coefficients of irreducible
corepresentations of $\kh_\sigma$. Since $\beta (B)\subset B\otimes_{alg}\kh_\sigma$, for any $b\in
B$ we may use the notation $\beta (b)=\sum_{uij}b^u_{ij}\otimes u_{ij}$ (finite
sum). See e.g. \cite{wo}. From the
coassociativity of $\beta$ we get
$$\sum_{uij}\beta
(b_{ij}^u)\otimes u_{ij}=\sum_{uijk}b_{ij}^u\otimes u_{kj}\otimes
u_{ik}$$
so $\beta (b_{ik}^u)=\sum_j b_{ij}^u\otimes u_{kj}$ for any
$u,i,k$. Thus
$$(id\otimes S)\beta (b^u_{ij})=(id\otimes S)(\sum_s b_{is}^u\otimes
u_{js})=\sum_s  b_{is}^u\otimes u_{sj}^*$$
Also $\widehat{\alpha}(1\otimes
u_{ij})=\sum_k1\otimes u_{ik}\otimes u_{kj}$, so we get that
$(\beta\odot\widehat{\alpha})(\beta (b)_{13})$ is equal to
$$\sum_{uij}(\sum_k 1\otimes 1\otimes u_{ik}\otimes
u_{kj})(\sum_s b_{is}^u\otimes 1\otimes 1\otimes
u_{sj}^*)=\sum_{uijks}b_{is}^u\otimes 1\otimes u_{ik}\otimes
u_{kj}u_{sj}^*$$
By summing over $j$ the last term is replaced by $(uu^*)_{ks}=\delta_{k,s}1$. Thus
$$(\beta\odot\widehat{\alpha})(\beta
(b)_{13})=\sum_{uik}b_{ik}^u\otimes 1\otimes u_{ik}\otimes 1=(\beta
(b)_{13})\otimes 1$$

From (a,b,c) we get that $(B\otimes (A\rtimes_\alpha\dh_\sigma
))^{\beta\odot\widehat{\alpha}}$, which is the fixed point algebra of $\sp\{ (1\otimes
1\otimes\kh_\sigma )\cdot \beta (B)_{13}\cdot \alpha (A)_{23}\}$ under
  $\beta\odot\widehat{\alpha}$, is $\sp\{ \beta  (B)_{13}\cdot \alpha
  (A)_{23}\}$. This finishes the proof of the first assertion and proves the non-degeneracy of the diagram in the statement. For the commuting square condition, remark that this diagram is the dual of the square on the left in the following diagram
$$\begin{matrix}
B &\subset & (B\otimes P)^{\beta\odot\pi} &\subset & B\otimes P\cr 
\cup &\ &\cup &\ &\cup\cr 
{\c} &\subset & P^\pi &\subset & P
\end{matrix}$$
where $P=A\rtimes_{\alpha}\dh_\sigma$ and $\pi =\widehat{\alpha}$. Since $\pi$ is dual, lemma 4.1 in
\cite{kac} applies and shows that the square on the right is a
non-degenerate commuting square. We also know that the
rectangle is a non-degenerate commuting square. Thus if we denote by $E_X:B\otimes P\to
B\otimes P$ the conditional expectation onto $X$, for any $X$, then
for any $b\in B$ we have $E_{P^\pi}(b)=E_P(b)=E_{{\c}}(b)$. This
proves the commuting square condition.
\end{proof}

We recall from \cite{aut} that any finite dimensional $\c^*$-algebra $B$ has a certain distingiushed trace, called canonical trace, which makes $\c\subset B$ a Markov inclusion. This is by definition the trace obtained by restricting the unique normalised trace of $\l (B)$ to the copy of $B$ given by the left regular representation. See \cite{aut} for other equivalent definitions (e.g. in terms of weights) and interpretations of this trace.

We denote by $\alg$ the category having as objects the finite dimensional $\c^*$-algebras and having as arrows the inclusions of $\c^*$-algebras which preserve the canonical traces. It is easy to see from definitions that any arrow in $\alg$ is Markov (see e.g. \S 6 in \cite{kac}). If $\kl$ is a Kac algebra we denote by $\kl -\alg$ the category having as objects the pairs $(A,\alpha )$ where $A$ is a finite dimensional $\c^*$-algebra and $\alpha :A\to A\otimes\kl$ is a coaction leaving invariant the canonical trace. The arrows in $\kl -\alg$ are by definition the inclusions in $\alg$ which commute with the coactions of $A$.

If $\kh$ is a compact Kac algebra and $(B,\beta )\in \kh-\alg$ and $(A,\alpha )\in \dh_\sigma -\alg$ then the $\c^*$-algebra $(B\otimes (A\rtimes_\alpha\dh_\sigma ))^{\beta\odot\widehat{\alpha}}$ in proposition 2.1 is finite dimensional, i.e. it is an object in $\alg$. We denote it by $B\Box_\kh A$. Note that the commuting square in proposition 2.1 shows that the inclusion $\c\subset B\Box_\kh A$ is Markov, i.e. that the trace on $B\Box_\kh A$ given by proposition 2.1 is its canonical trace.

If $(A^\prime ,\alpha^\prime )\subset (A,\alpha )$ is an arrow in $\dh_\sigma -\alg$ and $(B^\prime ,\beta^\prime )\subset (B,\beta )$ is an arrow in $\kh -\alg$ then proposition 2.1 gives a canonical embedding $B^\prime \Box_\kh A^\prime \subset B\Box_\kh A$. Since both $B^\prime \Box_\kh A^\prime$ and $B\Box_\kh A$ are endowed with their canonical traces, this inclusion has to be Markov, i.e. it is an arrow in $\alg$.

Summing up, we have constructed a bifunctor $\Box_\kh :\kh -\alg\times\dh_\sigma -\alg\to \alg$.

\begin{theo}
For any two arrows $B_0\subset B_1$ in $\kh -\alg$ and $A_0\subset A_1$ in $\dh_\sigma -\alg$
$$\begin{matrix}
B_0\Box_\kh A_1 &\subset & B_1\Box_\kh A_1\cr 
\cup &\ &\cup\cr 
B_0\Box_\kh A_0 &\subset & B_1\Box_\kh A_0
\end{matrix}$$
is a non-degenerate commuting square of finite dimensional von Neumann algebras.
\end{theo}

\begin{proof}
Step I. For $A_0=B_0=\c$ this follows from proposition 2.1.

Step II. We prove it for $A_0=\c$. Let $A=A_1$ and consider the following diagram
$$\begin{matrix}
A &\subset & B_0\Box_\kh A &\subset & B_1\Box_\kh A\cr 
\cup &\ &\cup &\ &\cup\cr 
\c &\subset & B_0 &\subset & B_1
\end{matrix}$$
By Step I, the square on the left and the rectangle are non-degenerate commuting squares. We want to prove that the square on the right is a non-degenerate commuting square. The non-degeneracy condition follows from
$$B_1\Box_\kh A=sp\{ A\cdot B_1\}\subset sp\{ B_0\Box_\kh A\cdot B_1\}$$
Let $x\in B_0\Box_\kh A$ and write $x=\sum_ib_ia_i$ with $b_i\in B_0$ and $a_i\in A$. Then
$$E_{B_1}(x)=\sum_ib_iE_{B_1}(a_i)=\sum_ib_iE_{\c}(a_i)=\sum_ib_iE_{B_0}(a_i)=E_{B_0}(x)$$
This proves the commuting square condition.

Step III. A similar proof shows that the proposition holds in the case $B_0=\c$.

Step IV. General case. We will use many times the following diagram
$$\begin{matrix}
A_1 &\subset & B_0\Box_\kh A_1 &\subset & B_1\Box_\kh A_1\cr 
\cup &\ &\cup &\ &\cup\cr 
A_0 &\subset & B_0\Box_\kh A_0 &\subset & B_1\Box_\kh A_0\cr
\cup &\ &\cup &\ &\cup\cr 
\c &\subset & B_0 &\subset & B_1
\end{matrix}$$
in which all the rectangles and all the squares, except possibly for the square in the statement, are non-degenerate
commuting squares (cf. Steps I, II, III). The non-degeneracy condition follows from
$$B_1\Box_\kh A_1=sp\{ A_1\cdot B_1\}\subset sp\{ B_0\Box_\kh A_1\cdot B_1\Box_\kh A_0\}$$
Let $x\in B_0\Box_\kh A_1$ and write $x=\sum_ib_ia_i$ with $b_i\in B_0$ and $a_i\in A_1$. Then
$$E_{B_1\Box_\kh A_0}(x)=\sum_ib_iE_{B_1\Box_\kh A_0}(a_i)=\sum_ib_iE_{A_0}(a_i)=\sum_ib_iE_{B_0\Box_\kh A_0}(a_i)=E_{B_0\Box_\kh A_0}(x)$$
This proves the commuting square condition. 
\end{proof}

We show now that the bifunctor $\Box_\kh$ behaves well with respect to basic constructions. This should be related to A. Wassermann's Invariance Principle \cite{was}.

If $\kl$ is a Kac algebra, a sequence of two arrows $A_0\subset A_1\subset A_2$ in $\kl -\alg$ is called a basic construction if $A_0\subset A_1\subset A_2$ is a basic construction in $\alg$ and if its Jones projection $e\in A_2$ is a fixed by the coaction $A_2\to A_2\otimes \kl$. An infinite sequence of basic constructions in $\kl -\alg$ is called a Jones tower in $\kl -\alg$.

\begin{prop}
If $B_0\subset B_1\subset B_2\subset B_3\subset\cdots$ is a Jones tower in $\kh -\alg$ and $A_0\subset A_1\subset A_2\subset A_3\subset\cdots$ is a Jones tower in $\dh_\sigma -\alg$ then
$$\begin{matrix}
\vdots &\ & \vdots &\ & \vdots &\ & \ \cr 
\cup &\ &\cup &\ &\cup\cr 
B_0\Box_\kh A_2 &\subset & B_1\Box_\kh A_2 &\subset & B_2\Box_\kh A_2 &\subset & \cdots\cr 
\cup &\ &\cup &\ &\cup\cr 
B_0\Box_\kh A_1 &\subset & B_1\Box_\kh A_1 &\subset & B_2\Box_\kh A_1 &\subset & \cdots\cr
\cup &\ &\cup &\ &\cup\cr 
B_0\Box_\kh A_0 &\subset & B_1\Box_\kh A_0 &\subset & B_2\Box_\kh A_0 &\subset & \cdots
\end{matrix}$$
is a lattice of basic constructions for non-degenerate commuting squares.
\end{prop}

\begin{proof}
We prove only that the rows are Jones towers (same proof for columns). By restricting attention to a pair of
consecutive inclusions, it is enough to prove the following statement:
if $B_0\subset B_1\subset B_2$ is a basic construction in $\kh -\alg$ and
$A$ is an object of $\dh_\sigma -\alg$ then $B_0\Box_\kh A\subset B_1\Box_\kh
A\subset B_2\Box_\kh A$ is a basic construction in $\alg$. For, we will use many times the following diagram
$$\begin{matrix}
A &\subset & B_0\Box_\kh A &\subset & B_1\Box_\kh A &\subset & B_2\Box_\kh A\cr 
\cup &\ &\cup &\ &\cup &\ &\cup\cr 
\c &\subset & B_0 &\subset & B_1 &\subset & B_2
\end{matrix}$$
in which all squares and rectangles are non-degenerate commuting
squares (cf. theorem 2.1). We use the abstract characterisation of the basic construction: $N\subset M\subset P$ is a basic construction with Jones projection $e\in P$ if and only if: (1) $P=sp\{ M\cdot e\cdot M\}$; (2) $[e,N]=0$; (3) $exe=E_N(x)e$ for any $x\in M$ and (4) $tr(xe)=\lambda tr(x)$ for any $x\in M$, where $\lambda$ is the inverse of the index
of $N\subset M$. Let $e\in B_2$ be the Jones projection for the basic construction $B_0\subset B_1\subset B_2$. With
$N=B_0$, $M=B_1$ and $P=B_2$ the verification of (1--4) is as follows.

(1) follows from the following computation:
$$B_2\Box_\kh A=sp\{ B_2\cdot A\} =sp\{ B_1\cdot e\cdot B_1\cdot A\} =sp\{ B_1\cdot e\cdot
B_1\Box_\kh A\}$$

(2) follows from $B_0\Box_\kh A=sp\{ B_0\cdot A\}$, from $[e,A]=0$ and from $[e,B_0]=0$.

(3) Let $x\in B_1\Box_\kh A$ and write $x=\sum_ib_ia_i$ with
$b_i\in B_1$ and $a_i\in A$. Then
$$exe=\sum_ieb_ia_ie=\sum_ieb_iea_i=\sum_iE_{B_0}(b_i)ea_i=\sum_iE_{B_0}(b_i)a_ie$$
$$E_{B_0\Box_\kh A}(x)e=\sum_iE_{B_0\Box_\kh A}(b_ia_i)e=\sum_iE_{B_0\Box_\kh A}(b_i)a_ie=\sum_iE_{B_0}(b_i)a_ie$$

(4) With the above notations, we have that
$$E_{B_2}(xe)=\sum_iE_{B_2}(b_ia_ie)=\sum_ib_iE_{B_2}(a_i)e=\sum_ib_iE_{\c}(a_i)e$$
and that $b_iE_{\c}(a_i)\in B_1$ for every $i$, so
$$tr_{B_2\Box_\kh A}(xe)=tr_{B_2}(E_{B_2}(xe))=\lambda\sum_itr_{B_1}(b_iE_{\c}(a_i))$$
$$tr_{B_1\Box_\kh A}(x)=tr_{B_1}(E_{B_1}(x))=\sum_itr_{B_1}(b_iE_{B_1}(a_i))=\sum_itr_{B_1}(b_iE_{\c}(a_i))$$
and this proves (4).
\end{proof}

\section{Commuting squares containing $\c$}

If $P$ and $Q$ are two unital $*$-algebras and $u\in P\otimes Q$ is a unitary we define a morphism of unital $*$-algebras $\iota_u:P\to P\otimes Q$ by $p\mapsto u(p\otimes 1)u^*$.

If $\kh$ is a compact Kac algebra we denote by $H\subset\kh$ the Hopf $*$-algebra consisting of coefficients of finite dimensional unitary corepresentations. We can associate unitary corepresentations and coactions of $\dh$ to representations of $H$ as follows. 
Let $V\in \b (l^2(\kh )\otimes l^2(\kh ))$ be the multiplicative unitary associated to $\kh$ and let $V^\prime\in M(\bar{H}\otimes \widehat{\bar{H}}_{red})$ be its version constructed in corollary A.6 in \cite{bs}, where $\bar{H}$
is the enveloping ${\c}^*$-algebra of $H$. Then any finite dimensional $*$-representation $\pi :H\to M_k$
has a unique continuous extension $\pi :\bar{H}\to M_k$, so one can define a corepresentation $\check{\pi}$ of
$\dh$ by $\check{\pi}=(\pi\otimes id )V^\prime$ and a coaction $\iota_{\check{\pi}} :M_k\to M_k\otimes \dh$ by $x\mapsto \check{\pi}(x\otimes 1)\check{\pi}^*$. See \cite{bs}.

\begin{lemm}
Let $v\in M_n\otimes H$ be a unitary corepresentation and let $\pi :H_\sigma\to M_k$ be a $*$-representation. Consider the objects $(M_n, \iota_v)\in \kh -\alg$ and $(M_k,\iota_{\check{\pi}})\in \dh_\sigma -\alg$ and form the corresponding algebra $M_n\Box_\kh M_k$. Then there exists an isomorphism
$$\begin{pmatrix}
M_k &\subset & M_n\Box_\kh M_k\cr
\cup &\ &\cup\cr
\c &\subset & M_n\end{pmatrix}\simeq
\begin{pmatrix}
{\c}\otimes M_k &\subset & M_n\otimes M_k\cr
\cup &\ &\cup\cr
{\c} &\subset &u(M_n\otimes{\c})u^*\end{pmatrix}$$
sending $z\mapsto 1\otimes z$ for any $z\in M_k$ and $y\mapsto\iota_u (y)$ for any $y\in M_n$, where $u=(id\otimes\pi )v$.
\end{lemm}

\begin{proof}
Consider the $*$-morphism
$$\Phi :M_n\otimes M_k\to M_n\otimes M_k\otimes\b (l^2(\kh_\sigma )),\,\,\,\,\, x\mapsto ad(v_{13}\check{\pi}_{23}u_{12}^*)(x\otimes 1)$$
Since both squares in the statement are non-degenerate commuting squares, all the assertions are consequences of the formulae $\Phi (1\otimes z)=\iota_{\check{\pi}} (z)_{23}$ for any $z\in M_k$ and $\Phi (\iota_u (y))=\iota_v(y)_{13}$ for any $y\in M_n$, to be proved now. The second one follows from the following computation
$$\Phi (u(y\otimes 1)u^*)=v_{13}(y\otimes 1\otimes 1)v_{13}^*=(v(y\otimes 1)v^*)_{13}$$
For the first formula, what we have to prove is that
$$v_{13}\check{\pi}_{23}u_{12}^*(1\otimes z\otimes 1)u_{12}\check{\pi}_{23}^*v_{13}^*=(\check{\pi}(z\otimes 1)\check{\pi}^*)_{23}$$
By moving the unitaries to the left and to the right we have to prove that
$$\check{\pi}_{23}^*v_{13}\check{\pi}_{23}u_{12}^*\in (\c\otimes M_k\otimes\c )^\prime =M_n\otimes\c\otimes \b (l^2(\kh_\sigma ))$$
Call this unitary $U$. Since $\check{\pi}=(\pi\otimes id )V^\prime$ we have
$$U=(id\otimes\pi\otimes id)(V_{23}^*v_{13}V_{23}v_{12}^*)$$
By \cite{bs} the comultiplication of $H_\sigma$ is given by $\Delta :y\mapsto V^*(1\otimes y)V$. On the other hand since $v$ is a
corepresentation of $H$, its adjoint $v^*$ is a corepresentation of $H_\sigma$, i.e. we have $(id\otimes\Delta )(v^*)=v^*_{12}v^*_{13}$. We get that
$$V_{23}^*v_{13}V_{23}=(V_{23}^*v_{13}^*V_{23})^*=((id\otimes\Delta )(v^*))^*=(v^*_{12}v^*_{13})^*=v_{13}v_{12}$$
Thus $U=v_{13}$ and we are done.
\end{proof}

If $V,W$ are $\c$-linear spaces, with $V$ finite dimensional, the canonical isomorphism $\l (V,V\otimes W)\simeq  \l (V)\otimes W$ is denoted by $\beta\mapsto u_\beta$. Its inverse is denoted by $u\mapsto \beta_u$.

If $S$ is a commuting square we denote by $S^\uparrow$ and $S^\to$ the commuting squares obtained by performing an upwards (resp. to the right) basic construction to $S$.

We denote by $\Box_u$ the vertex model associated to a biunitary $u$.

We recall from \cite{ver} that if $H$ is a Hopf $*$-algebra, $v\in M_n\otimes H$ is a unitary corepresentation and $\pi :H\to M_k$ is a $*$-representation then the element $(id\otimes\pi )v\in M_n\otimes M_k$ is a biunitary. Conversely, given a biunitary $u\in M_n\otimes M_k$, the category of triples $(H,v,\pi )$ such that $(id\otimes\pi )v=u$ has a universal object, called the minimal model for $u$.

For $D\in\alg$ we make use of the canonical embedding $D\subset \l (D)$.

\begin{lemm}
Let $A,B\in\alg$ and $T$ be a commuting square of the form
$$\begin{matrix}
\c\otimes\l (A)&\subset & B\otimes \l (A)\cr
\cup &\ &\cup \varphi\cr
\c &\subset &B\end{matrix}$$

(i) The element $u_\varphi\in\l (B)\otimes \l (A)$ is a biunitary.

(ii) If $(H,v,\pi )$ is the minimal model for $u_\varphi$ then $\beta_v$ is a coaction.

(iii) We have $\varphi (b)=\iota_{u_\varphi}(b)$ in $\l (B)\otimes \l (A)$, for any $b\in B$.

(iv) $T^\to$ is isomorphic to $\Box_{u_\varphi}$.
\end{lemm}

\begin{proof}
(i) Denote by $\mu_B:B\otimes B\to B$ the multiplication and by $\eta_B:\c\to B$ the unit. We recall that the universal $*$-algebra $A^{aut}(B)$ generated by the entries of a unitary matrix $V$ such that $\eta_B\in Hom(1,V)$ and $\mu_B\in Hom(V^{\otimes 2},V)$ has a unique structure of Hopf $*$-algebra which makes $V\in \l (B)\otimes A^{aut}(B)$ a corepresentation of it. The element $\beta_V$ is a coaction. See \cite{wang2}, \cite{aut}.

Since $\varphi$ is a $*$-morphism we get that $\eta_B\in Hom(1,u_\varphi)$ and that $\mu_B\in Hom(u_\varphi^{\otimes 2},u_\varphi)$ with the notations in \cite{aut} (same proof as for lemma 1.2 (i,ii) in there). Since $T$ is a commuting square we have that $(tr\otimes id)\varphi =tr (.)1$, so the same computation as in the first part of the proof of lemma 1.2 (iv) in \cite{aut} shows that $u^*_\varphi u_\varphi =1$, i.e. that $u_\varphi$ is a unitary. It follows that there exists a representation $\nu :A^{aut}(B)\to \l (A)$ such that $(id\otimes\nu )V=u_\varphi$. Thus $u_\varphi$ is biunitary and $(A^{aut}(B),V,\nu )$ is a model for it (see \cite{ver}).

(ii) By using the universal property of the minimal model we get a morphism $p:A^{aut}(B)\to H$ such that $(id\otimes p)V=v$. It follows that $\beta_v$ is equal to $(id\otimes p)\beta_V$, which is a coaction.

(iii) Since $v$ is a corepresentation we have $\beta_v(b)=\iota_v(b)$ in $\l (B)\otimes H$, for any $b\in B$. By applying $id\otimes\pi$ we get that $(id\otimes\pi )\beta_v(b)=(id\otimes\pi )\iota_v(b)$ in $\l (B)\otimes \l (A)$, for any $b\in B$. The term on the right is $\iota_{(id\otimes\pi )v}(b)=\iota_{u_\varphi}(b)$ and the term on the left is $\beta_{(id\otimes\pi )v}(b)=\beta_{u_\varphi}(b)=\varphi (b)$.

(iv) Consider the following diagram
$$\begin{matrix}
\l (A) &\subset & B\Box_{\kh}\l (A) &\subset & \l (B)\Box_{\kh}\l (A)\cr 
\cup &\ &\cup &\ &\cup\cr 
{\c} &\subset & B &\subset & \l (B)
\end{matrix}$$
associated to the objects $(B,\beta_v)$ and $(\l (B),\iota_v)$ of $\kh -\alg$ and $(\l (A),\iota_{\check{\pi}} )$ of $\dh_\sigma -\alg$. By proposition 2.2 this is a basic construction for non-degenerate commuting squares. By applying lemma 3.1 with $(H,v,\pi )$ we get an isomorphism between the big square in the above diagram and the big square in the diagram below
$$\begin{matrix}
\c\otimes\l (A) &\subset & B\otimes\l (A) &\subset & \l (B)\otimes \l (A)\cr 
\cup &\ &\cup &\ &\cup\cr 
{\c} &\subset & \varphi (B) &\subset & u_\varphi (\l (B)\otimes\c )u_\varphi^*
\end{matrix}$$
which sends $z\mapsto 1\otimes z$ for any $z\in \l (A)$ and $y\mapsto \iota_{u_\varphi}(y)$ for any $y\in \l (B)$. It follows that this diagram is also a basic construction for non-degenerate commuting squares. Thus $T^\to =\Box_{u_\varphi}$.
\end{proof}

We recall from \cite{js} that any commuting square consisting of $\c$ and three algebras of matrices is isomorphic to a vertex model. Thus if $S$ is a non-degenerate commuting square of the form
$$\begin{matrix}
A &\subset & X\cr
\cup &\ &\cup\cr
\c &\subset & B\end{matrix}$$
then $S^{\uparrow\rightarrow}$ is isomorphic to a vertex model. The choice of a biunitary $u\in\l (B)\otimes \l (A)$ such that $S^{\uparrow\rightarrow}\simeq\Box_u$ is not unique. We will need the following (still not unique) special choice for such a biunitary.

(1) The commuting square $S^\uparrow$ consists of the algebras $\c$, $B$, $\l (A)$ and $B\otimes \l (A)$. Since any embedding $\l (A)\subset B\otimes \l (A)$ is conjugate to the canonical embedding given by $m\mapsto 1\otimes m$, we get an isomorphism of the form
$$S^\uparrow\simeq T:=\begin{pmatrix}
\c\otimes\l (A)&\subset & B\otimes \l (A)\cr
\cup &\ &\cup \varphi\cr
\c &\subset &B\end{pmatrix}$$
where $\varphi :B\subset B\otimes \l (A)$ is a certain inclusion. By performing a basic construction to the right this induces an isomorphism $S^{\uparrow\rightarrow}\simeq T^\to$. On the other hand, by applying lemma 3.2 we get an isomorphism $T^\to\simeq \Box_{u_\varphi}$. By composing these two isomorphisms we get an isomorphism $S^{\uparrow\rightarrow}\simeq \Box_{u_\varphi}$.

(2) Consider the dual commuting square $\widehat{S}$, obtained from $S$ by interchanging $A$ and $B$. By applying (1) to $\widehat{S}$ we get an inclusion $\phi :A\subset A\otimes \l (B)$ and an isomorphism $\widehat{S}^{\uparrow\to}\simeq \Box_{u_\phi}$. Since the dual of $\widehat{S}^{\uparrow\to}$ is  $S^{\uparrow\rightarrow}$ (obvious) and the dual of $\Box_U$ is $\Box_{\sigma (U)^*}$ for any biunitary $U$, where $\sigma$ is the flip (see the beginning of appendix I below) we get by dualising an isomorphism $S^{\uparrow\rightarrow}\simeq \Box_{\sigma (u_\phi )^*}$.

(3) By composing the isomorphisms in (1) and (2) we get an isomorphism $\Box_{u_\varphi}\simeq \Box_{\sigma (u_\phi )^*}$. It follows that $u_\varphi =ad(r\otimes q)(\sigma (u_\phi )^*)$ for some unitaries $r\in \l (B)$ and $q\in \l (A)$. We define $u\in\l (B)\otimes \l (A)$ to be $ad(r\otimes 1)(\sigma (u_\phi )^*)=ad(1\otimes q^*)(u_\varphi )$.

For unexplained terminology in the statement (ii) below see appendix II below.

\begin{theo}
Let $S$ be a non-degenerate commuting square of the form
$$\begin{matrix}
A &\subset & X\cr
\cup &\ &\cup\cr
\c &\subset & B\end{matrix}$$
Let $u\in\l (B)\otimes \l (A)$ be as above and let $(H,v,\pi )$ be its minimal model. 

(i) The maps $\beta_v :B\to B\otimes\kh$ and $\beta_{\check{\pi}}:A\to A\otimes\dh_\sigma$ are coactions and $S$ is isomorphic to the commuting square
$$\begin{matrix}
A &\subset & B\Box_\kh A\cr
\cup &\ &\cup\cr
\c &\subset & B\end{matrix}$$
associated to the objects $(B,\beta_v )\in\kh -\alg$ and $(A, \beta_{\check{\pi}})\in\dh_\sigma -\alg$.

(ii) The product type coaction $\gamma_\pi$ associated to $\pi$ is outer, we have $Z(B)\cap B^{\beta_v}=\c$ and the vertical subfactor associated to $S$ is isomorphic to $R\subset (B\otimes (R\rtimes_{\gamma_\pi} \dh_\sigma ))^{\beta_v\odot \widehat{\gamma_\pi}}$.
\end{theo}

\begin{proof}
(i) Step I. Since $(H,v,\pi )$ is the minimal model for $u$, $(H,v, ad(q)\pi)$ is the minimal model for $ad(1\otimes q)u=u_\varphi$. Thus lemma 3.2 shows that $\beta_v$ is a coaction.

Step II. We prove that $\beta_{\check{\pi}}$ is a coaction. We use terminology and notations from the appendix I below. By using the construction $\pi\mapsto\check{\pi}$ in the beginning of this section we can define an embedding $i:H^\circ_\sigma\subset \dh_\sigma$ by $f\pi \mapsto (f\otimes id)\check{\pi}$ for any coefficient $f\pi$. We have $(id\otimes i)\widehat{\pi}=\check{\pi}$ for any representation $\pi$. Let us apply now the result in Step I with the dual of $S$ instead of $S$ and with $\sigma (u)^*$ as the choice of the biunitary. We get that if $(K,w,\rho )$ is the minimal model for $\sigma (u)^*$ then $\beta_w$ is a coaction of $K$. By applying theorem 4.1 below we get an embedding $\phi :K\subset H^\circ_\sigma$ such that $(id\otimes\phi )w=\widehat{\pi}$. It follows that $\beta_{\widehat{\pi}}$ is a coaction of $H^\circ_\sigma$, so $\beta_{\check{\pi}} =(id\otimes i)\beta_{\widehat{\pi}}$ is a coaction of $\dh_\sigma$.

Step III. It remains to prove the isomorphism part. For, consider the diagram
$$\begin{matrix}
\l (A) &\subset & B\Box_\kh \l (A) &\subset & \l (B)\Box_\kh\l (A)\cr 
\cup &\ &\cup &\ &\cup\cr 
A &\subset & B\Box_\kh A &\subset & \l (B)\Box_\kh A\cr
\cup &\ &\cup &\ &\cup\cr 
\c &\subset & B &\subset & \l (B)
\end{matrix}$$
associated to the objects $(B,\beta_v)$ and $(\l (B),\iota_v)$ of $\kh -\alg$ and $(A,\beta_{\check{\pi}})$ and $(\l (A),\iota_{\check{\pi}} )$ of $\dh_\sigma -\alg$. By proposition 2.2 this is a lattice of basic construction for non-degenerate commuting squares. On the other hand, by applying lemma 3.1 with $(H,v,\pi )$ we get an isomorphism from the big square in this diagram, say $Q$, to $\Box_u$, which sends $z\mapsto 1\otimes z$ for any $z\in \l (A)$ and $y\mapsto \iota_u(y)$ for any $y\in \l (B)$. By composing it with the isomorphism $\Box_u\simeq S^{\uparrow\rightarrow}$ we get an isomorphism $\Phi :Q\simeq S^{\uparrow\to}$. Since $\Phi (A)=A$ and $\Phi (B)=B$, $\Phi$ implements an isomorphism between the lower left square in $Q$ and $S$.

(ii) Proposition 2.2 shows that inclusion in (ii) is the vertical subfactor associated to the commuting square in (i). Let us prove that $\gamma_\pi$ is outer. For, we apply theorem 5.1. Since $\pi$ is inner faithful (cf. \cite{ver}) it remains to prove that $\pi$ contains the counit. For, we come back to Lemma 3.2. With the notations in there, we have that $\eta_B\in Hom(1,V)$, so in particular $V$ contains a copy of the trivial corepresentation $1$. Since $v =(id\otimes p)V$, we have that $1\in v$. By arguing as in Step III in the above proof of (i), with the notations in there, we have $1\in w$. Since $\widehat{\pi}=(id\otimes\phi )w$, it follows that $1\in\widehat{\pi}$ and by dualising we get that $\pi$ contains the counit. Finally, the ergodicity condition $Z(B)\cap B^{\beta_v}=\c$ follows from factoriality and from theorem 4.1 in \cite{aut}. 
\end{proof}

In \cite{kac} the standard invariant of any fixed point subfactor of the form $P^\pi\subset (B\otimes P)^{\beta\odot\pi}$ is computed and shown to depend only on the coaction $\beta$. We call it the Popa system of $\beta$. By applying theorem 5.1 in \cite{kac} to the subfactor in theorem 4.1 (ii) we get:

\begin{coro}
The standard invariant of the vertical subfactor associated to $S$ is isomorphic to the Popa system of the coaction $\beta_v$.
\end{coro}

This generalises the computation in \cite{ver} for vertex models.

\section{Appendix I: Duality for vertex models}

The theory of commuting squares has a duality given by
$$\begin{pmatrix}
A&\subset & X\cr
\cup &\ &\cup\cr
Y &\subset &B\end{pmatrix}\longleftrightarrow\begin{pmatrix}
B&\subset & X\cr
\cup &\ &\cup\cr
Y &\subset &A\end{pmatrix}$$

At the level of vertex models, the canonical isomorphism
$$\begin{pmatrix}
{\c}\otimes M_k &\subset & M_n\otimes M_k\cr
\cup &\ &\cup\cr
{\c} &\subset &u(M_n\otimes{\c})u^*
\end{pmatrix}\simeq 
\begin{pmatrix}
\sigma (u)^*(M_k\otimes {\c})\sigma (u) &\subset & M_k\otimes M_n\cr
\cup &\ &\cup\cr
{\c} &\subset &{\c}\otimes M_n
\end{pmatrix}$$
where $\sigma$ is the flip, shows that this duality is implemented by the duality
$$M_n\otimes M_k \ni u\longleftrightarrow\sigma (u)^*\in M_k\otimes M_n$$
for biunitaries. This gives rise to a duality for minimal models for biunitaries. Here we find an abstract formulation for it. We use freely terminology and notations from \cite{ver}.

\begin{defi}
A Hopf algebra $H$ is said to be subreflexive if there exists
an embedding $i_H:H\subset H^{\circ\circ}$ making the following diagram commutative:
$$\begin{matrix}
H &\subset & H^{\circ\circ}\cr 
\cap &\ &\cap\cr 
H^{**} &\to & H^{\circ *}
\end{matrix}$$

A duality between two subreflexive Hopf algebras $H$ and $K$ is a
pair of embeddings $\phi :K\subset H^{\circ}_\sigma$ and $\psi
:H\subset K^{\circ}_\sigma$ such that $\psi =\phi^{\circ}i_H$ and $\phi =\psi^{\circ}i_K$.
\end{defi}

Any finite dimensional Hopf algebra $H$ is subreflexive and $H^*_\sigma$ is its unique dual. By using theorem 4.1 below and \cite{ver} we get that for $u_1,\ldots ,u_n\in {\bf U}(k)$ the Hopf algebras ${\c}[<u_1,\ldots ,u_n>]$ and ${\mathcal R}_\c (\overline{<u_1,\ldots ,u_n>})$ are subreflexive and dual to each other. This example shows that the dual is not unique in general.

If $v\in M_n\otimes H$ is a corepresentation, $v=\sum e_{kl}\otimes v_{kl}$, then the dual representation $\widehat{v}:H^{\circ}\to M_n$ is $\widehat{v}:\psi\mapsto\sum \psi (v_{kl})e_{kl}$. If $\pi :H\to M_k$ is a representation, $\pi :x\mapsto \sum \psi_{ij}(x)f_{ij}$, then the dual corepresentation $\widehat{\pi}\in M_k\otimes H^{\circ}$ is $\widehat{\pi}=\sum f_{ij}\otimes \psi_{ij}$.

\begin{theo}
Let $u$ be a biunitary and let $(H,v,\pi )$ and $(K,w,\rho )$ be the minimal models for $u$ and for $\sigma (u )^*$. Then $H$ and $K$ are subreflexive and there exists a unique duality of subreflexive Hopf
algebras $\phi :K\subset H^{\circ}_\sigma$ and $\psi
:H\subset K^{\circ}_\sigma$ such that $(id\otimes \phi )w=\widehat{\pi}$,
$(id\otimes \psi )v=\widehat{\rho}$, $\widehat{v}\phi =\rho$ and
$\widehat{w}\psi =\pi$.
\end{theo}

\begin{proof}
We prove first that $H$ is subreflexive (same proof for $K$). We have
$$(id\otimes \widehat{v})\widehat{\pi}=\sum
\psi_{ij}(v_{kl})f_{ij}\otimes e_{kl}=\sigma (\sum
\psi_{ij}(v_{kl})e_{kl}\otimes f_{ij})=\sigma ((id\otimes\pi )v)$$
Thus if $u=(id\otimes\pi )v$ then $(H^{\circ},
\widehat{\pi},\widehat{v})$ is a model for $\sigma (u)$. By using this
twice we get that $(H^{\circ\circ},\widehat{\widehat{v}}, 
\widehat{\widehat{\pi}})$ is a model for $u$. Since $(H,v,\pi )$ is the minimal model for $u$, we get a morphism $p:H\to H^{\circ\circ}$
satisfying $(id\otimes p)v=\widehat{\widehat{v}}$. Consider now the diagram
$$\begin{matrix}
H &\to^p & H^{\circ\circ}\cr 
\cap &\ &\cap\cr 
H^{**} &\to & H^{\circ *}
\end{matrix}$$
The condition $(id\otimes p)v=\widehat{\widehat{v}}$ shows that this
diagram commutes on the coefficients of $v$. Since these coefficients
generate $H$ as a Hopf algebra, this diagram commutes. It remains to prove that $p$ is injective. By using the above diagram
it is enough to prove that the composition $H\subset H^{**}\to H^{\circ *}$
is injective. The kernel of this map is by definition $(H^{\circ})^{\perp
  (H)}$. Since the coefficients of $\pi$ generate a dense subalgebra of $H^*$ (see \cite{ver}), we have in
particular that $H^{\circ}$ is dense in $H^*$ and this gives $(H^{\circ})^{\perp
  (H)}=\{ 0\}$.

We construct now the duality. Since $(H^{\circ},\widehat{\pi},\widehat{v})$ is a model for $\sigma (u)$, the triple $(H^{\circ}_\sigma ,
\widehat{\pi}^*,\widehat{v})$ is a model for $\sigma (u)^*$. By definition of $\widehat{v}$ its space of coefficients ${\mathcal
  C}_{\widehat{v}}\subset H^{\circ\circ}$ is the image of the space ${\mathcal
  C}_v\subset H$ by the embedding $H\subset H^{\circ\circ}$. On the other hand
  for any subreflexive Hopf algebra $H$ we have from definitions that
  $H^{\perp (H^{\circ})}=\{ 0\}$, so $H$ is dense in $H^{\circ *}$. Together with
  the fact that ${\mathcal
  C}_v$ generates $H$ as a Hopf algebra, we get that ${\mathcal
  C}_{\widehat{v}}$ generates in $H^{\circ\circ}$ a Hopf subalgebra which is
  dense in $H^{\circ *}$. Thus if we define ${\mathcal C}_\pi\subset H^{\circ}_\sigma$ to be the
space of coefficients of $\pi$ and denote by $<{\mathcal
  C}_\pi>\subset H^{\circ}_\sigma$ the Hopf subalgebra generated by it, then $(<{\mathcal
  C}_\pi>, \widehat{\pi}^*,\widehat{v}_{\mid <{\mathcal
  C}_\pi>})$ is a bi-faithful model for $\sigma (u)^*$. By \cite{ver} we get that this triple is isomorphic to the minimal model $(K,w,\rho )$. We get in this way an isomorphism
$K\simeq <{\mathcal
  C}_\pi>$ and by composing it with the embedding $<{\mathcal
  C}_\pi>\subset H^{\circ}_\sigma$ we get an embedding $\phi :K\subset H^{\circ}_\sigma$. The
same method gives also an embedding $\psi
:H\subset K^{\circ}_\sigma$ and all the assertions are clear with these $\phi$ and
$\psi$.
\end{proof}

\section{Appendix II: Outerness of product type coactions}

If $P$ is a finite von Neumann algebra, a coaction $\pi
:P\to P\otimes\kh$ is said to be minimal if $(P^\pi )^\prime\cap
P={\c}$ and if it faithful. If $R$ is a finite von Neumann algebra, a coaction
$\gamma :R\to R\otimes\dh$ is said to be outer if $R^\prime\cap
(R\rtimes_\gamma\dh )={\c}$, i.e. if its dual coaction
$\widehat{\gamma} :R\rtimes_\gamma\dh\to (R\rtimes_\gamma\dh
)\otimes\kh$ is minimal. See e.g \cite{kac}.

If $u\in M_n(\kh )$ is a unitary corepresentation and $\pi
:P\to P\otimes\kh$ is a coaction then $\pi_u: x\mapsto
u_{13}((id\otimes\pi )x)u_{13}^*$ is a coaction of $\kh$ on
$M_n\otimes P$. Note that for the trivial coaction $\iota :\c\to\c\otimes\kh$ we have $\iota_u (x)= u(x\otimes 1)u^*$. See e.g. proposition 2.1 in \cite{kac}.

\begin{lemm}
Let $P$ be a finite von Neumann algebra and $\pi :P\to P\otimes\kh$ be
a faithful coaction. Then $\pi$ is minimal if and only if the relative commutant of $P^\pi\subset
(M_n\otimes P)^{\pi_w}$ is isomorphic to $End(w)$, for any unitary
corepresentation $w\in M_n(H)$.
\end{lemm}

\begin{proof}
Note first that the equality $M_n^{\iota_w}=End(w)$ gives an embedding
of $End(w)$ in the relative commutant of in the statement. The ``only
if'' part is clear from theorem 5.1 in \cite{kac} (because $P^\pi\subset (M_n\otimes
P)^{\pi_w}$ is a Wassermann-type subfactor) and will not be used
here. For the converse, let $Irr(\kh )$ be a complete system of non-equivalent unitary irreducible corepresentations of $\kh$, each of them co-acting on some ${\c}^k$ (see e.g. \cite{wo}). For any $p\in {\mathcal P}$ we
write $\pi (p)=\sum_{uij}p^u_{ij}\otimes u_{ij}$ and we use the
formula $\pi (p_{ij}^u)=\sum_k p_{kj}^u\otimes u_{ki}$, which follows
from the coassociativity of $\pi$. For any $u\in Irr(\kh )$ we define its
spectral projection $E_u:P\to P$ by $p\mapsto dim(u)(id\otimes h)(\pi
(p)(1\otimes \chi (u)^*))$. With the above notation, we have $E_u(p)=\sum_ip_{ii}^u$, so it
follows that $\{ E_u\}_{u\in Irr({\mathcal A} )}$ are orthogonal projections
with respect to the trace of $P$, that their images $P^u=Im(E_u)$ are in
${\mathcal P}$ and that ${\mathcal P}$ decomposes as a direct sum $\oplus_u P^u$ (see \cite{kac} for details).

Step I. We prove that an element $p\in P$ commutes with
$P^{\pi}$ if and only if all its spectral projections commute with
$P^\pi$. The ``if'' part is clear. Conversely, for any $x\in
P^\pi$ we have $\pi (xp)=\alpha (px)$, so $\pi (p)$ commutes
with $x\otimes 1$. Thus
$$xE_u(p)=dim(u)(id\otimes h)((x\otimes 1)\pi (p)(1\otimes\chi
(u)^*))=E_u(p)x$$
for any unitary corepresentation $u\in Irr(\kh )$.

Step II. Assume that $\pi$ is not minimal and choose a non-scalar element
$q\in (P^\pi )^\prime\cap P$. By Step I all its spectral projections
commute with $P^\pi$. The condition in the statement gives for
$m=1$ and $w=1$ that $P^\pi$ is a factor, so the spectral
projection $E_1(q)$ is a scalar. Thus there exists a non-trivial
irreducible corepresentation, say $u\in M_n(\kh )$ and a
nonzero element $p\in (P^\pi )^\prime\cap P^u$. If we denote $\pi (p)=\sum p_{ij}\otimes u_{ij}$ then we have $\pi
(p_{ij})=\sum_kp_{kj}\otimes u_{ki}$. This shows that the matrix
$m=(p_{ji})_{ij}$ is a (nonzero) $u$-eigenmatrix, i.e. that $m$ satisfies
$(id\otimes\pi )m=m_{12}u_{13}$. It follows that
$$\begin{pmatrix} {1}&{0}\cr {0}&{u_{13}}\end{pmatrix}
\begin{pmatrix} {0}&{(id\otimes\pi )m}\cr {0}&{0}\end{pmatrix}
\begin{pmatrix} {1}&{0}\cr {0}&{u^*_{13}}\end{pmatrix}
=\begin{pmatrix} {0}&{m\otimes 1}\cr {0}&{0}\end{pmatrix}$$
Consider now the unitary corepresentation
$$u^+:=(n\otimes 1)\oplus u=\begin{pmatrix} 1&0\cr 0&u\end{pmatrix}
\in M_2(M_n\otimes \kh )=M_2\otimes M_n\otimes\kh$$
The above equality shows that the matrix $\begin{pmatrix}
  0&m\cr 0&0\end{pmatrix}$ is in the fixed point algebra
$X_u=(M_2\otimes M_n\otimes P)^{\pi_{u^+}}$ (note that this is the
algebra appearing in the quantum extension of Wassermann's extension
of Connes' $2\times 2$ matrix trick, lemma 7.2 in \cite{kac}). Thus
$$\begin{pmatrix} 0& m\cr 0&0\end{pmatrix}\in (P^\pi
)^\prime\cap X_u=End(u^+)\otimes 1=\begin{pmatrix} M_n&0\cr 0&{\c}
  I\end{pmatrix}$$
where we use the hypothesis with $w=u^+$. Thus $m=0$, contradiction.
\end{proof}

To any representation $\pi :H\to M_k$ one can associate a product type coaction $\gamma_\pi$ of $\dh$ on the hyperfinite $II_1$-factor $R$
in the following way. Consider the corepresentation $\check{\pi}\in
M_k\otimes\dh$ (see \S 3) and for any $n$ form the tensor product of
corepresentations
$\check{\pi}_n=\check{\pi}\otimes\overline{\check{\pi}}\otimes\check{\pi}\otimes\ldots$,
where $\overline{\check{\pi}}$ is the complex conjugate of
$\check{\pi}$. Then $\iota_{\check{\pi}_{n+1}}$ extends
$\iota_{\check{\pi}_n}$ for any $n$, so we may consider the inductive
limit coaction $\iota_{\check{\pi}_\infty}$ on the inductive limit
$M_k^{\otimes\infty}$. By extending $\iota_{\check{\pi}_\infty}$ to the adherence $R=\overline{M_k^{\otimes\infty}}^w$ we
get the coaction $\gamma_\pi$.

We say that a representation $\pi :H\to M_k$ contains the counit if the counit $\varepsilon\in H^\circ$ of $H$ is a coefficient of
$\pi$ and that it is inner faithful if its space of coefficients is dense in $H^*$ (see \S 1 in \cite{ver} for terminology and notations). Note that any finite dimensional Kac algebra has a finite dimensional inner faithful representation containing the counit (take the regular one) so the result below generalises \cite{y}.

\begin{theo}
The product type coaction associated to an inner faithful representation containing the counit is outer.
\end{theo}

\begin{proof}
Let $\pi :H\to M_k$ be our representation and $\gamma_\pi :R\to R\otimes\dh$ be its product type coaction. By applying lemma 5.1 we have to prove the following statement: for any corepresentation $w\in M_n(H)$ the relative
commutant of the inclusion $R\subset (M_n\otimes (R\rtimes_{\gamma_\pi}\dh ))^{\widehat{\gamma_\pi}_w}$ is $End(w)$. We may and will assume that the coefficients of $w$ generate $H$. It is easy to see that this inclusion is isomorphic to $R\subset (M_n\otimes (R\rtimes_{\gamma_\pi}\dh_\sigma ))^{\iota_{w^t}\odot\widehat{\gamma_\pi}}$, where $w^t\in M_n\otimes\kh_\sigma$ is the corepresentation equal to the transpose of $w$ and where we use notations from \S 1 (see \cite{kac} for details). By using lemma 3.1 and proposition 2.2, this inclusion is isomorphic to the vertical subfactor associated to the vertex model corresponding to the biunitary $u=(id\otimes\pi )w^t$. Theorem 5.1 in \cite{ver} applies and shows that the relative commutant is $End(y)$, where $(G,y,\nu )$ is the minimal model for
$$u_{12}u^{*t}_{13}=((id\otimes\pi )w^t)_{12}((id\otimes t\pi )\overline{w})_{13}=(id\otimes\rho )w^t$$
where $\rho =(\pi\otimes t\pi S)\sigma\Delta$. Since $\pi$ contains the
counit, the spaces of coefficients of both $\pi$ and $t\pi S$
are contained in the space of coefficients of $\rho$. It follows that
$\rho$ is inner faithful, so $(H_\sigma ,w^t,\rho )$ is the minimal model for
$u_{12}u^{*t}_{13}$. Thus $(G,y,\nu )=(H_\sigma ,w^t,\rho )$, so our
commutant is $End(w^t)=End(w)$.
\end{proof}

\end{document}